\documentclass[12pt]{amsart}

\usepackage{graphicx}
\usepackage{amsmath}

\def\R{\mathbb{R}}

\def\cD{\mathcal{D}}

\def\al{\alpha}

\def\ga{\gamma}

\def\la{\lambda}

\def\si{\sigma}

\def\om{\omega}

\def\Si{\Sigma}

\def\Om{\Omega}

\newcommand{\der}{{\rm d}}

\numberwithin{equation}{section}

\newtheorem{theorem}{Theorem}[section]

\newtheorem{proposition}[theorem]{Proposition}

\theoremstyle{remark}

\theoremstyle{remark}

\usepackage{amssymb,stmaryrd}
\usepackage{amscd}

\newcommand{\qr}{
\begin{center}
\includegraphics[scale=0.5]{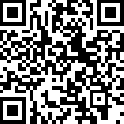}\\

\text{Scan the QR code to view more articles from the author}
\end{center}}

\author{Matthew Randall}
\address{Institute of Mathematical Sciences \\
ShanghaiTech University\\
393 Middle Huaxia Road\\
Shanghai 201210\\
China}
\email{mjrandall@shanghaitech.edu.cn}

\title{A Monge normal form for the rolling distribution}


\begin{document}

\begin{abstract}
Using a parametrisation of $sl_2$ given by the second prolongation of the group action of unimodular fractional linear transformations as presented in an article of Clarkson and Olver \cite{co96}, we find a Monge normal form describing the rolling of two hyperboloid surfaces over each other. 
\end{abstract}

\maketitle

\pagestyle{myheadings}
\markboth{Randall}{A Monge normal form for the rolling distribution}

\section{Introduction}

Let $\cD$ be a maximally non-integrable rank 2 distribution on a 5-manifold $M$. The maximally non-integrable condition of $\cD$ determines a filtration of the tangent bundle $TM$ given by
\[
\cD \subset [\cD,\cD] \subset [\cD,[\cD,\cD]]\cong TM.
\]
The distribution $[\cD, \cD]$ has rank 3 while the full tangent space $TM$ has rank 5, hence such a geometry is also known as a $(2,3,5)$-distribution. Let $M_{xyzpq}$ denote the 5-dimensional mixed order jet space $J^{2,0}(\R,\R^2) \cong J^2(\R,\R)\times \R$ with local coordinates given by $(x,y,z,p,q)=(x,y,z,y',y'')$ (see also \cite{tw13}, \cite{tw14}). Let $\cD_{F(x,y,z,y',y'')}$ denote the maximally non-integrable rank 2 distribution on $M_{xyzpq}$ associated to the underdetermined differential equation $z'=F(x,y,z,y',y'')$. This means that the distribution is annihilated by the following three 1-forms
\begin{align*}
\om_1=\der y-p \der x, \qquad \om_2=\der p-q \der x, \qquad  \om_3=\der z-F(x,y,z,p,q) \der x.
\end{align*}
Such a distribution $\cD_{F(x,y,z,y',y'')}$ is said to be in Monge normal form (see page 90 of \cite{tw13}). In Section 5 of \cite{conf}, it is shown how to associate canonically to such a $(2,3,5)$-distribution a conformal class of metrics of split signature $(2,3)$ (henceforth known as Nurowski's conformal structure or Nurowski's conformal metrics) such that the rank 2 distribution is isotropic with respect to any metric in the conformal class. The method of equivalence \cite{cartan1910} (also see the introduction to \cite{annur}, Section 5 of \cite{conf} and \cite{Strazzullo}) produces the 1-forms $(\theta_1, \theta_2,\theta_3, \theta_4, \theta_5)$ that gives a coframing for Nurowski's metric. These 1-forms satisfy the structure equations
\begin{align}\label{cse}
\der \theta_1&=\theta_1\wedge (2\Om_1+\Om_4)+\theta_2\wedge \Om_2+\theta_3 \wedge \theta_4,\nonumber\\
\der \theta_2&=\theta_1\wedge\Om_3+\theta_2\wedge (\Om_1+2\Om_4)+\theta_3 \wedge \theta_5,\nonumber\\
\der \theta_3&=\theta_1\wedge\Om_5+\theta_2\wedge\Om_6+\theta_3\wedge (\Om_1+\Om_4)+\theta_4 \wedge \theta_5,\\
\der \theta_4&=\theta_1\wedge\Om_7+\frac{4}{3}\theta_3\wedge\Om_6+\theta_4\wedge \Om_1+\theta_5 \wedge \Om_2,\nonumber\\
\der \theta_5&=\theta_2\wedge \Om_7-\frac{4}{3}\theta_3\wedge \Om_5+\theta_4\wedge\Om_3+\theta_5\wedge \Om_4,\nonumber
\end{align}
where $(\Om_1, \ldots, \Om_7)$ and two additional 1-forms $(\Om_8, \Om_9)$ together define a rank 14 principal bundle over the 5-manifold $M$ (see \cite{cartan1910} and Section 5 of \cite{conf}). A representative metric in Nurowski's conformal class \cite{conf} is given by
\begin{align}\label{metric}
g=2 \theta_1 \theta_5-2\theta_2 \theta_4+\frac{4}{3}\theta_3 \theta_3.
\end{align}
When $g$ has vanishing Weyl tensor, the distribution is called maximally symmetric and has split $G_2$ as its group of local symmetries. For further details, see the introduction to \cite{annur} and Section 5 of \cite{conf}. For further discussion on the relationship between maximally symmetric $(2,3,5)$-distributions and the automorphism group of the split octonions, see Section 2 of \cite{tw13}.

$(2,3,5)$-distributions also arise from the study of the configration space of two surfaces rolling without slipping or twisting over each other \cite{AN14}, \cite{BH} and \cite{BM}. The configuration space can be realised as the An-Nurowski circle twistor distribution \cite{AN14} and in the case of two spheres with radii in the ratio $1:3$ rolling without slipping or twisting over each other, there is again maximal $G_2$ symmetry. 

In the work of \cite{r17}, a description of maximally symmetric $(2,3,5)$-distributions obtained from Pfaffian systems with $SU(2)$ symmetry was discussed and its relationship with the rolling distribution was investigated. In particular, the An-Nurowski circle twistor bundle can be realised by considering the Riemannian surface element of the unit sphere arising from one copy of $SU(2)$ and the other Riemannian surface element with Gaussian curvature $9$ or $\frac{1}{9}$ from another copy of $SU(2)$. Both Lie algebras of $su(2)$ are parametrised by the left-invariant vector fields. See \cite{r17} for further details. 

In the $SL(2)$ case, the spheres are replaced by hyperboloids equipped with Lorentzian signature metrics, but the ratios of the Gaussian curvatures of the surfaces remain unchanged for there to be maximal split $G_2$ symmetry. It is in this context where the Lorentzian surface elements of the hyperboloids arise from two different copies of $SL(2)$ that we consider the rolling distribution. We base the computations here from that already obtained in \cite{r17}. 

Using a different parametrisation of $sl_2$ given in \cite{co96} obtained from the second prolongation of the group action of unimodular fractional linear transformations, we show a change of coordinates that bring the 1-forms annihilating the rolling distribution to the 1-forms encoding the Monge equation
\begin{align}\label{meq}
\om_1&=\der y-p \der x,\nonumber\\
\om_2&=\der p-q \der x,\nonumber\\
\om_3&=\der z-\left(qz^2+\frac{1}{\al^2-1}(\sqrt{q}z-\frac{1}{2\sqrt{q}x})^2\right) \der x. 
\end{align} 
In Monge normal form, $F(x,y,z,p,q)=\left(qz^2+\frac{1}{\al^2-1}(\sqrt{q}z-\frac{1}{2\sqrt{q}x})^2\right)$. This is the content of Theorem \ref{mainthm} in the paper.
The distribution has maximal split $G_2$ symmetry whenever $\al^2=9$ or $\al^2=\frac{1}{9}$.
The Monge normal form we obtain here is of course not unique and is defined up to symmetries of the distribution $\cD_{F}$. Some more work would be needed to see if we can find a further change of coordinates to eliminate the $z^2$ term in the expression of $F(x,y,z,p,q)$. 
The computations here are done with the aid of the \texttt{DifferentialGeometry} package in MAPLE 2018. 

\section{A hyperboloid rolling distribution}
The real $SL(2)$ analogue of Section 2 in \cite{r17} is given as follows. 
We can parametrise the Lie algebra $sl_2$ using the left-invariant vector fields, from which we get
\begin{align*}
\si_1&=\sinh(y) \cosh(z) \der p-\sinh(z) \der y,\\
\si_2&=-\sinh(y) \sinh(z) \der p+\cosh(z) \der y,\\
\si_3&=-\der z-\cosh(y) \der p.
\end{align*}
These 1-forms satisfy the relations
\begin{align*}
\der \si_1&=-\si_2\wedge \si_3,\\
\der \si_2&=-\si_1\wedge \si_3,\\
\der \si_3&=\si_1\wedge \si_2,
\end{align*}
and the Lorentzian hyperbolic surface element with Gauss curvature $-1$ is given by
\[
\si_2\si_2-\si_1\si_1=\der y^2-\sinh^2(y)\der p^2.
\]
We find that the distribution annihilated by the the Pfaffian system spanned by the 1-forms
\begin{align}\label{rolling}
\om_1=-(\si_1+\exp(\al x)\der q),\quad \om_2=\si_2+\der x,\quad \om_3=-(\si_3+\al \exp(\al x)\der q)
\end{align}
is a $(2,3,5)$-distribution whenever $\al^2 \neq 1$ and has maximal symmetry whenever $\al^2=\frac{1}{9}$ or $\al^2=9$. 
The 1-forms can be completed into a coframing by taking $\om_4=-\der x$ and $\om_5=\exp(\al x)\der q$. To satisfy Cartan's structure equations, we take
\begin{align*}
\theta_1&=\om_1,\quad \theta_2=\om_2,\quad \theta_3=K^{\frac{1}{3}}\om_3,\\
\theta_4&=K^{-\frac{1}{3}}\om_4+P\theta_1+Q\theta_2+R\theta_3,\\
\theta_5&=K^{-\frac{1}{3}}\om_5+S\theta_1+T\theta_2+U\theta_3,
\end{align*}
where $K=\frac{1}{\al^2-1}$, $P=T=0$, $R=U=0$ and $Q=S=\frac{3 \al^2-7}{10(\al^2-1)^{\frac{2}{3}}}$.
This gives the conformal metric
\begin{align*}
K^{\frac{1}{3}}g&=2\om_1\om_5-2\om_2\om_4+\frac{3\al^2-7}{5}K(\om_1\om_1-\om_2\om_2)+\frac{4}{3}K\om_3\om_3\\
&=2 \om_4^2-2\om_5^2+\frac{1}{2}(\si_1-\om_5)^2-\frac{1}{2}\om_1^2-\frac{1}{2}(\si_2+\om_4)^2+\frac{1}{2}\om_2^2\\
&\hspace{12pt}+\frac{3\al^2-7}{5(\al^2-1)}(\om_1^2-\om_2^2)+\frac{4}{3(\al^2-1)}\om_3^2.
\end{align*}
Let
\begin{align*}
\bar \om_1=-\si_1+\om_5\quad \mbox{and} \quad \bar\om_2=\si_2+\om_4. 
\end{align*}
Using the fact that
\begin{align*}
\si_1^2-\si_2^2+\om_5^2-\om_4^2=\frac{1}{2}\left(\om_1^2-\om_2^2\right)+\frac{1}{2}(\bar \om_1^2-\bar \om_2^2),
\end{align*}
we find
\begin{align}\label{diagonalm}
K^{\frac{1}{3}}g&=2 \om_4^2-2\om_5^2+\frac{1}{2}\bar \om_1^2-\frac{1}{2}\om_1^2-\frac{1}{2}\bar \om_2^2+\frac{1}{2}\om_2^2+\frac{3\al^2-7}{5(\al^2-1)}(\om_1^2-\om_2^2)+\frac{4}{3(\al^2-1)}\om_3^2\nonumber\\
&=2 \om_4^2-2\om_5^2+\frac{1}{2}\bar \om_1^2+\frac{1}{2}\om_1^2-\frac{1}{2}\bar \om_2^2-\frac{1}{2}\om_2^2+(\frac{3\al^2-7}{5(\al^2-1)}-1)(\om_1^2-\om_2^2)+\frac{4}{3(\al^2-1)}\om_3^2\nonumber\\
&=\om_4^2-\om_5^2-(\si_2^2-\si_1^2)+(\frac{3\al^2-7}{5(\al^2-1)}-1)(\om_1^2-\om_2^2)+\frac{4}{3(\al^2-1)}\om_3^2.
\end{align}
In the terminology of \cite{r17}, the 1-forms $\{\bar \om_1,\bar \om_2, \om_3\}$ form a ``sign-reversed" Pfaffian system corresponding to the symmetry $(\al,x,q) \mapsto (-\al,-x,-q)$. 
Observe that the metric we obtain in (\ref{diagonalm}) is diagonal. Each hyperboloid surface $\Si$ and $\Si'$ is equipped with a Lorentzian surface element. The Gauss curvature of the surface element $\om_4^2-\om_5^2$ over $\Si$ is given by $-\al^2$, and the Gauss curvature of $\si_2^2-\si_1^2$ over $\Si'$ is given by $-1$. The maximally symmetric $(2,3,5)$-distribution is obtained when the ratios of the Gauss curvature of the two surfaces are $9$ or $\frac{1}{9}$. In these cases, the metric (\ref{diagonalm}) has vanishing Weyl tensor. See also \cite{AN14}.

We call this distribution {\it rolling} because the 1-forms $\{\om_1, \om_2,\om_3\}$ given in (\ref{rolling}) define a connection on the An-Nurowski circle twistor bundle over the product of the two surfaces $\Si$ and $\Si'$ with elements (or metric tensor) given by $\om_4^2-\om_5^2$ and $\si_2^2-\si_1^2$ respectively. Furthermore, the An-Nurowski circle twistor bundle realises the system of two surfaces rolling without slipping or twisting over each other \cite{AN14}. See \cite{r17} for details in the $SU(2)$ case. We shall also refer to the Pfaffian system spanned by the 1-forms in (\ref{rolling}) as a {\it rolling} system. 

We now reparametrise the Lie algebra of $sl_2$ using the vector fields arising from the second prolongation of the group action of unimodular fractional linear transformations. From this reparametrisation, we are able to reduce the distribution to a Monge normal form after a change of coordinates. 

\section{A Monge normal form for the rolling distribution}
Let
\begin{align*}
\tau_1=\der y+y \der z,\quad \tau_2=-(\der p-p \der z),\quad \tau_3=-\der z.
\end{align*}
The set of 1-forms
\begin{align*}
s_1=\tau_1+y^2 \tau_2,\quad s_2=\tau_2,\quad s_3=\tau_3-2 y \tau_2,
\end{align*}
forms a basis dual to the $sl_2$ Lie algebra of vector fields
\begin{align*}
\partial_y, \quad y^2\partial_y-2y \partial_z-(2 y p+1)\partial_p, \quad y\partial_y-\partial_z-p\partial_p, 
\end{align*}
(see \cite{co96}) that arise from the second prolongation of the group action of unimodular fractional linear transformations.
We have
\begin{align*}
\der s_1&=-s_1\wedge s_3,\\
\der s_2&=s_2\wedge s_3,\\
\der s_3&=-2s_1\wedge s_2.
\end{align*}
We now consider the {\it rolling} distribution annihilated by the 1-forms
\begin{align*}
\theta_1&=s_1+s_2-\exp(\al x)\der q=\der y+y \der z-(y^2+1)(\der p-p \der z)-\exp(\al x)\der q,\\
\theta_2&=s_1-s_2+\der x=(1-y^2)(\der p-p \der z)+(\der y+y \der z)+\der x,\\
\theta_3&=-s_3+\al \exp(\al x)\der q=-2 y(\der p-p \der z)+dz+\al \exp(\al x)\der q.
\end{align*}
The Pfaffian system is equivalent to the Pfaffian system (\ref{rolling}) above
\begin{align*}
\om_1=-(\si_1+\exp(\al x)\der q),\quad \om_2=\si_2+\der x,\quad \om_3=-(\si_3+\al \exp(\al x)\der q)
\end{align*}
through the isomorphism
\begin{align*}
\{\si_1,\si_2,\si_3\}\cong \{-(s_1+s_2),s_1-s_2,-s_3\}.
\end{align*}
The distribution is equivalently annihilated by the 1-forms
\begin{align}
\theta_1+\frac{1}{\al}\theta_3&=\der y+y \der z-(y^2+1+\frac{2}{\al}y)(\der p-p \der z)+\frac{1}{\al}\der z,\nonumber\\
\theta_1-\theta_2&=-2(\der p-p \der z)-\der x-\exp(\al x)\der q,\label{new1forms}\\
\theta_3-y(\theta_1-\theta_2)&=\der z+y \der x+(y+\al)\exp(\al x)\der q.\nonumber
\end{align}
From the equations (\ref{new1forms}) determined by the ideal $\{\theta_1,\theta_2,\theta_3\}$, we are now going to find the change of coordinates that bring it to Monge normal form. 
We define
\begin{align*}
\hat y&=\exp(z)(y+\frac{1}{\al}),\quad \hat p=\exp(-z) p,\quad \hat q=\frac{1}{\al}\exp(-\al x),\\
\hat z&=\exp(z-\al x),\quad \hat x=q-\frac{1}{\al}\exp(-\al x).
\end{align*}
Under this change of coordinates, the system spanned by the 1-forms given in (\ref{new1forms}) is equivalently spanned by the 1-forms
\begin{align*}
\der \hat y-\left(\exp(2 z)(y^2+\frac{2 y}{\al}+1)\right)\der \hat p&=\der \hat y-\left(\hat y^2+(1-\frac{1}{\al^2})\exp(2 z)\right) \der \hat p\\
&=\der \hat y-\left(\hat y^2+(1-\frac{1}{\al^2})\frac{\hat z^2}{\al^2 \hat q^2}\right) \der \hat p,\\
\der\hat p+\frac{1}{2 \hat z}\der \hat x\quad \mbox{and} \quad \der \hat z+(\hat y+\frac{\left(1-\frac{1}{\al^2}\right)\hat z}{\hat q})\der \hat x.
\end{align*}
We used the fact that $\exp(2 z)=\frac{\hat z^2}{\al^2 \hat q^2}$.
Now consider further the map
\[
(x,y,z,p,q)=\left(2 \hat z, \hat x+2\hat z \hat p,\hat y,\hat p,\frac{1}{4(\hat z \hat y+(1-\frac{1}{\al^2})\frac{\hat z^2}{\hat q})}\right).
\]
We find
\begin{align*}
\der y-p \der x&=\der \hat x+2 \hat p \der \hat z+2 \hat z \der \hat p-2\hat p\der \hat z=\der \hat x+2 \hat z \der \hat p,\\
\der p-q \der x&=\der \hat p-\frac{2}{4(\hat z \hat y+(1-\frac{1}{\al^2})\frac{\hat z^2}{\hat q})}\der \hat z=-\frac{1}{2\hat z}(\der \hat x+\frac{1}{\hat y+(1-\frac{1}{\al^2})\frac{\hat z}{\hat q}}\der \hat z),\\
\der z-F \der x&=\der \hat y-2F\der \hat z=\der \hat y-4 F(\hat z\hat y+(1-\frac{1}{\al^2})\frac{\hat z^2}{\hat q})\der \hat p,
\end{align*}
where $F=F(x,y,z,p,q)$ is to be determined. It follows that 
\[
F= q z^2+(1-\frac{1}{\al^2})q \exp(2z).
\]
Using the inverse map
\[
(\hat x,\hat y,\hat z,\hat p,\hat q)=\left(y-xp, z,\frac{x}{2}, p,\left(1-\frac{1}{\al^2}\right)\frac{q x^2}{1-2 z q x}\right),
\]
we deduce that
\begin{align*}
\frac{\hat z}{\hat q}=\frac{x}{2\hat q}=\left(1-\frac{1}{\al^2}\right)^{-1}\left(\frac{1}{2qx}-z\right).
\end{align*}
This gives
\begin{align*}
\exp(2z)=\frac{x^2}{4\al^2\hat q^2}=\frac{1}{\al^2}\left(1-\frac{1}{\al^2}\right)^{-2}\left(\frac{1}{2qx}-z\right)^2.
\end{align*}
We therefore obtain
\begin{align*}
F&=qz^2+\frac{1}{\al^2-1}(\sqrt{q}z-\frac{1}{2\sqrt{q}x})^2\\
&=q z^2+\frac{1}{\al^2-1}q\left(z^2-\frac{z}{q x}+\frac{1}{4 q^2 x^2}\right).
\end{align*}
When $\al^2=9$, we obtain
\begin{align}\label{max1a}
F=\frac{9}{8}q z^2-\frac{1}{8}\frac{z}{x}+\frac{1}{32q x^2}.
\end{align}
When $\al^2=\frac{1}{9}$, we obtain
\begin{align}\label{max1b}
F=-\frac{1}{8}q z^2+\frac{9}{8}\frac{z}{x}-\frac{9}{32q x^2}.
\end{align}
In both of these cases they give the Monge normal form of a maximally symmetric $(2,3,5)$-distribution. 
We have the following theorem
\begin{theorem}\label{mainthm}
Let $\cD$ be the $(2,3,5)$-distribution associated to the {\it rolling} system spanned by the 1-forms $\{\om_1, \om_2, \om_3\}$ given by (\ref{rolling}). Then by the change of coordinates given above, this Pfaffian system can be brought into the Monge normal form given by
 \begin{align*}
\om_1&=\der y-p \der x,\\
\om_2&=\der p-q \der x,\\
\om_3&=\der z-\left(qz^2+\frac{1}{\al^2-1}(\sqrt{q}z-\frac{1}{2\sqrt{q}x})^2\right) \der x.
\end{align*} 
The $(2,3,5)$-distribution has split $G_2$ symmetry whenever $\al=\pm\frac{1}{3}$, $\al=\pm3$.
\end{theorem}

\section{$SL(2)$ Pfaffian systems}
In this Section, we give an example of a homogeneous bracket-generating $(2,3,5)$-distribution that can be reduced to the above Monge normal form given by the {\it rolling} distribution. This $(2,3,5)$-distribution can be seen as a generalisation of the {\it rolling} distribution, where we assume that each hyperboloid surface comes from a copy of $SL(2)$. This is also the real analogue of the $SU(2)$ picture discussed in Section 5 of \cite{r17}. However, because our parametrisation of $sl_2$ is different and does not use the left-invariant vector fields on $SL(2)$, the geometric relationship with the {\it rolling} distribution is not immediately apparent. 

On $M^5$ with local coordinates given by $(x,y,z,p,q)$, consider the 1-forms
\begin{align*}
\om_1&=\der y+y \der z,\quad \om_2=-(\der p-p \der z),\quad \om_3=-\der z,\\
\om_4&=\der q+q \der z,\quad \om_5=-(\der x-x \der z).
\end{align*}
The vector fields dual to the set of 1-forms
\begin{align*}
s_1=\om_1+y^2 \om_2,\quad s_2=\om_2,\quad s_3=\om_3-2 y \om_2,
\end{align*}
form a copy of the Lie algebra $sl_2$, and
the vector fields dual to the set of 1-forms
\begin{align*}
s_4=\om_4+q^2 \om_5,\quad s_5&=\om_5,\quad s_6=\om_3-2 q \om_5,
\end{align*}
form a second copy of $sl_2$. We form the following set of 1-forms
\begin{align*}
s_1&=\om_1+y^2 \om_2,\quad s_2=\om_2,\quad s_3=\om_3-2 y \om_2,\\
\bar{s}_3&=\om_3-2 y \om_2-2 q\om_5=s_3-2 q \om_5,\quad \bar{s}_4=\frac{1}{q}(\om_4+q^2\om_5),\quad \bar{s}_5=q \om_5.
\end{align*}
These 1-forms satisfy the equations
\begin{align*}
\der s_1&=- s_1\wedge \bar s_3-2 s_1 \wedge \bar s_5,\\
\der s_2&=s_2\wedge \bar s_3+2 s_2 \wedge \bar s_5,\\
\der \bar s_3&=-2s_1\wedge s_2-2 \bar s_4 \wedge \bar s_5,\\
\der \bar s_4&=\bar s_4\wedge \bar s_5,\\
\der \bar s_5&=\bar s_4 \wedge \bar s_5. 
\end{align*}
Let 
\begin{align}
\theta_1&=s_1-\beta \bar s_4,\nonumber\\
\theta_2&=s_2-\ga \bar s_5,\label{gen1forms}\\
\theta_3&=s_3+\bar s_4+\bar s_5.\nonumber
\end{align}
\begin{proposition}
The distribution given by the kernel of the 1-forms $\{\theta_1,\theta_2,\theta_3\}$ in (\ref{gen1forms}) is a bracket-generating $(2,3,5)$-distribution whenever $\beta\ga \neq 1$.
\end{proposition}

We shall call the Pfaffian system given by the 1-forms $\{\theta_1,\theta_2,\theta_3\}$ an $SL(2)$ Pfaffian system. It can be seen that these three 1-forms can be be completed into a coframing that satisfies Cartan's structure equations. 
Let
\begin{align*}
\theta_1&=s_1-\beta \bar s_4,\\
\theta_2&=s_2-\ga \bar s_5,\\
\bar \theta_3&=K^\frac{1}{3}(s_3+\bar s_4+\bar s_5),\\
\theta_4&=K^{-\frac{1}{3}}(\beta \bar s_4),\\
\theta_5&=K^{-\frac{1}{3}}(-\ga \bar s_5)+T \theta_2,
\end{align*}
where 
\begin{align*}
K=\frac{\beta \ga}{2(\beta \ga-1)} \mbox{~and~} T=\frac{-9\beta \ga\pm3\pm\sqrt{3(3\beta^2\ga^2-26\beta \ga+3)}}{6(\beta \ga-1)(K)^{\frac{1}{3}}}
\end{align*}
are constants. Then Cartan's structure equations (\ref{cse}) are satisfied for the 1-forms $(\theta_1, \theta_2, \bar \theta_3, \theta_4, \theta_5)$.
We find that Nurowski's conformal metric given by
\begin{align*}
g=2 \theta _1\theta_5-2 \theta_2\theta_4+\frac{4}{3}\bar \theta_3\bar \theta_3
\end{align*}
is conformally flat whenever $\beta\ga=9$ or $\beta\ga=\frac{1}{9}$.
We have
\begin{align*}
K^{\frac{1}{3}}g=-2 \ga(1+\la) s_1\bar s_5-2 \beta(1+\la) s_2\bar s_4+2\beta\ga(2+\la) s_4s_5+2\la s_1s_2+\frac{4}{3}K(s_3+\bar s_4+\bar s_5)^2
\end{align*}
where
\[
\la=K^{\frac{1}{3}}T=\frac{-9\beta \ga\pm3\pm\sqrt{3(3\beta^2\ga^2-26\beta \ga+3)}}{6(\beta \ga-1)}.
\]
Also note that $\bar s_4\bar s_5=s_4s_5$.

\begin{theorem}
The $(2,3,5)$-distribution given by the kernel of the 1-forms $\{\theta_1,\theta_2,\theta_3\}$ in (\ref{gen1forms})
has split $G_2$ symmetry whenever $\beta \ga=9$ or $\beta \ga=\frac{1}{9}$.
\end{theorem}
We now encode the Pfaffian system given by (\ref{gen1forms}) by a Monge equation equivalent to (\ref{meq}). 
The 1-forms annihilating the vector fields can be expressed as
\begin{align*}
\theta_1&=(\der y+y \der z)-y^2(\der p-p \der z)-\frac{\beta}{q}((\der q+q \der z)-q^2(\der x-x \der z)),\\
\theta_2&=(p \der z-\der p)+\ga q (\der x-x\der z),\\
\theta_3&=2 y(\der p-p \der z)+\frac{1}{q}\der q-2 q (\der x- x \der z).
\end{align*}
Equivalently, taking linear combinations, the vector fields are annihilated by the 1-forms
\begin{align*}
\theta_1+\beta\theta_3+\frac{\beta}{\ga}\theta_2&=(\der \tilde y+\tilde y \der z)-\left(\frac{\tilde y}{\ga}+(\tilde y+(\beta-\frac{1}{\ga}))(\tilde y-\beta)\right)(\der p-p \der z),\\
\theta_2&=(p \der z-\der p)+\ga q (\der x-x\der z),\\
\frac{1}{2 q ( \ga y-1)}\left(\theta_3+2 y \theta_2\right)&= (\der x- x \der z)+\frac{\der q}{2(\tilde y+(\beta-\frac{1}{\ga}))\ga q^2},
\end{align*}
where $\tilde y=y-\beta$. 
We we now take 
\begin{align*}
\hat y&=\frac{1}{c} \exp(z)\tilde y,\\
\hat p&=c \exp(-z) p,\\
\hat x&=-c \exp(-z) x,\\
\hat q&=-\frac{1}{\ga q},\\
\hat z&=(\beta-\frac{1}{\ga})\frac{\exp(z)}{c}.
\end{align*}
This reduces the system (\ref{gen1forms}) to the ideal spanned by the 1-forms
\begin{align*}
&\der \hat x-\hat q \der \hat p,\\
&\der\hat y-\left(\hat y^2-\beta\left(\beta-\frac{1}{\ga}\right)\frac{\exp(2z)}{c^2}\right)\der \hat p=\der \hat y-(\hat y^2-\frac{\beta}{\beta-\frac{1}{\ga}} \hat z^2) \der \hat p, \\
&\der \hat x-\frac{c \exp(-z)}{2(\tilde y+\beta-\frac{1}{\ga})}\der \hat q=\der \hat x-\frac{1}{2(\hat y+\frac{\exp(z)}{c}(\beta-\frac{1}{\ga}))}\der \hat q=\der \hat x-\frac{1}{2(\hat y+\hat z)}\der \hat q.
\end{align*}
We now map into more familiar coordinates on the mixed jet space by taking $(x, y, z, p, q)=(\hat q, \hat x-\hat p\hat q, -\hat y,-\hat p, -\frac{1}{2(\hat y+\hat z)\hat q})$, from which we obtain
\begin{align*}
&\der y-p \der x=\der \hat x-\hat q\der \hat p-\hat p\der \hat q +\hat p \der \hat q=\der \hat x-\hat q\der \hat p, \\
&\der p-q \der x=-\der \hat p+\frac{1}{2(\hat y+\hat z)\hat q}\der \hat q=-\frac{1}{\hat q}(\der \hat x-\frac{1}{2(\hat y+\hat z)}\der \hat q),\\
&\der z-F\der x=-\der \hat y-2(\hat y+\hat z)\hat q F\der \hat p=-\der \hat y+(\hat y^2-\mu \hat z^2)\der \hat p,
\end{align*}
with $\mu=\frac{\beta \ga}{\beta \ga-1}$ and $F=F(x,y,z,p,q)$ again to be determined.
This means 
\begin{align*}
F&=-\frac{1}{2(\hat y+\hat z)\hat q}\left(\hat y^2-\mu\hat z^2\right)\\
&=q(z^2-\mu(z-\frac{1}{2qx})^2)\\
&=\left(1-\frac{\beta \ga}{\beta \ga-1}\right)z^2 q-\left(\frac{\beta \ga}{\beta \ga-1}\right)\frac{1}{4 q x^2}+\left(\frac{\beta \ga}{\beta \ga-1}\right)\frac{z}{x}\\
&=z^2 q-\frac{\beta \ga}{\beta \ga-1}(z^2q+\frac{1}{4 qx^2}-\frac{z}{x})\\
&=z^2q+\frac{1}{\frac{1}{\beta \ga}-1}\left(\sqrt{q}z-\frac{1}{2 \sqrt{q} x}\right)^2. 
\end{align*}
This is equivalent to the form given in (\ref{meq}) after identifying $\beta \ga=\al^{-2}$. When $\beta\ga=\frac{1}{9}$, this gives (\ref{max1a}), and when $\beta\ga=9$, this gives (\ref{max1b}).

\begin{theorem}\label{thm2}
Let $\cD$ be the $(2,3,5)$-distribution associated to the $SL(2)$ Pfaffian system spanned by the 1-forms $\{\theta_1, \theta_2, \theta_3\}$ in (\ref{gen1forms}). Then by the change of coordinates given above, this Pfaffian system can be brought into the Monge normal form equivalent to the {\it rolling} distribution given by
 \begin{align*}
\om_1&=\der y-p \der x,\\
\om_2&=\der p-q \der x,\\
\om_3&=\der z-\left(qz^2+\frac{\beta \ga}{1-\beta \ga}(\sqrt{q}z-\frac{1}{2\sqrt{q}x})^2\right) \der x.
\end{align*} 
This has split $G_2$ symmetry whenever $\beta \ga=\frac{1}{9}$ or $9$.
\end{theorem}

\section{Conformal metric from Monge normal form}
In this section, for the sake of completeness, we  give the description of the Nurowski's conformal class associated to the above Monge normal form (\ref{meq}). 
The distribution $\cD_{F}$ where $F=qz^2+\frac{1}{\al^2-1}(qz^2-\frac{z}{x}+\frac{1}{4qx^2})=qz^2+\frac{1}{\al^2-1}(\sqrt{q}z-\frac{1}{2\sqrt{q}x})^2$ is annihilated by the three 1-forms
\begin{align*}
\om_1=\der y-p \der x,\quad \om_2=\der p-q \der x,\quad \om_3=\der z-\left(qz^2+\frac{1}{\al^2-1}(qz^2-\frac{z}{x}+\frac{1}{4qx^2})\right) \der x.
\end{align*} 
These three 1-forms are completed to a coframing on $M_{xyzpq}$ by the additional 1-forms
\begin{align*}
\om_4=\frac{1}{2 q^3 x^2 (\al^2-1)}\der q-\frac{4 \al^2 q^2 x^2 z^2-4\al^2 q x z+(3-2\al^2)}{4(\al^2-1)^2 x^3 q^2}\der x, \quad \om_5=-\der x.
\end{align*}
If we take
\begin{align*}
\theta_1&=\om_3-\frac{4 \al^2 q^2 x^2 z^2-1}{4 q^2 x^2 (\al^2-1)}\om_2, \quad \theta_2=\om_1,\quad \theta_3=K^{\frac{1}{3}}\om_2,\\
\theta_4&=K^{-\frac{1}{3}}\om_4+a_{41}\theta_1+a_{42}\theta_2+a_{43}\theta_3,\\
\theta_5&=K^{-\frac{1}{3}}\om_5+a_{51}\theta_1+a_{52}\theta_2+a_{53}\theta_3,
\end{align*}
where 
\begin{align*}
K&=\frac{1}{2 q^3 x^2 (\al^2-1)}, \quad a_{41}=0,\\
a_{42}&=\frac{2^{\frac{1}{3}}\al^2(\al^2-9)(4 qxz(qxz-1)+1)}{60x^{\frac{10}{3}}(\al^2-1)^{\frac{8}{3}}q^2}, \\
a_{43}&=-\frac{2^{\frac{2}{3}}(12\al^2q^2x^2z^2-8\al^2qxz-2\al^2+3)}{12x^{\frac{5}{3}}(\al^2-1)^{\frac{4}{3}}q},\\
a_{51}&=0,\quad a_{52}=\frac{2^{\frac{1}{3}}(2\al^2-3)}{5(\al^2-1)^{\frac{2}{3}}x^{\frac{1}{3}}},\quad a_{53}=-2^{\frac{2}{3}}qx^{\frac{4}{3}}(\al^2-1)^{\frac{2}{3}},
\end{align*}
then we find Cartan's structure equations are satisfied for $(\theta_1, \theta_2, \theta_3, \theta_4, \theta_5)$, and the coframing
gives a representative metric from Nurowski's conformal class with 
\begin{align*}
g=2 \theta_1 \theta_5-2\theta_2 \theta_4+\frac{4}{3}\theta_3\theta_3. 
\end{align*}
The metric is conformally flat when $\al^2=9$ or $\al^2=\frac{1}{9}$, in which case the distribution has maximal split $G_2$ symmetry. In comparison to the metric obtained in (\ref{diagonalm}), the metric here is no longer in diagonal form.

\qr 

\end{document}